# Convergence analysis of the Chebyshev-Legendre spectral method for a class of Fredholm fractional integro-differential equations

A. Yousefi, S. Javadi, E. Babolian, E. Moradi

Department of Computer Science, Faculty of Mathematical Sciences and Computer, Kharazmi University, Tehran, Iran

javadi@khu.ac.ir, asadyosefi@gmail.com, babolian@khu.ac.ir, eslam.moradi@gmail.com,

**Abstract**

In this paper, we propose and analyze a spectral Chebyshev-Legendre approximation for fractional order integro-differential equations of Fredholm type. The fractional derivative is described in the Caputo sense. Our proposed method is illustrated by considering some examples whose exact solutions are available. We prove that the error of the approximate solution decay exponentially in $L^2$-norm.

## 1. Introduction

Many phenomena in engineering, physics, chemistry, and the other sciences may be applied by models using mathematical tools from fractional calculus. The theory of derivatives and integrals of fractional order allow us to describe physical phenomena more accurately [1-2]. Furthermore most problems cannot be solved analytically, and hence finding good approximate solution, using numerical methods will be very helpful. Recently, several numerical methods have been given to solve fractional differential equations (FDEs) and fractional integro-differential equations (FIDEs). These methods include collocation method [3-4], variational iteration method [5], Adomian decomposition method [6], Homotopy perturbation method [5,7], fractional differential transform method [8-9], the reproducing kernel method [10], and wavelet method [11-13].

Spectral methods are an emerging area in the field of applied sciences and engineering. These methods provide a computational approach that has achieved substantial popularity over the last three decades. They have been applied successfully to numerical simulations of many problems in fractional calculus ([14-20]).

In this paper, we are concerned with numerical solutions of the following equation:

$$\sum_{i=0}^{n} a_i D^{(i)} y(t) = f(t) + \int_0^1 k(t,s) D^\alpha y(s)\, ds,$$
$$m - 1 < \alpha \leq m, m \in \mathbb{N}, \quad t \in [0,1], \tag{1}$$

subject to the initial values

$$y^{(i)}(0) = d_i, \quad i = 0,1,\ldots,n-1, \tag{2}$$

where $D^\alpha$ is the fractional derivative in the Caputo sense, $f(t)$ and $k(t,s)$ are the known functions that are supposed to be sufficiently smooth and $d_i$ for any $i$ is constant. Existence and uniqueness of the solution of the Eq. (1) have been shown in [21]. The authors in [22] applied the backward and central-difference formula for approximating solution at the mesh points.

The fractional derivative are global, i.e. they are defined over the whole interval $I = [0,1]$, and therefore global method, such as spectral methods, are better suited for FDEs and FIDEs. Yousefi and et al. [20] introduced a quadrature shifted Legendre tau method based on the Gauss-Lobatto interpolation for solving Eq. (1). Inspired by the work of [23-24], we extend the approach to Eq. (1) and provide a rigorous convergence analysis for the Chebyshev-Legendre method. We show that approximate solutions are convergent in $L^2$−norm.

The structure of this paper is as follows: In section 2, some necessary definitions and mathematical tools of the fractional calculus which are required for our subsequent developments are introduced. In section 3, the Chebyshev-Legendre method of FIDEs is obtained. The rest of this section is devoted to apply the proposed method for solving Eq. (1) by using the shifted Legendre and Chebyshev polynomials. After this section, we discuss about convergence analysis and then, some numerical experiments are presented in Section 5 to show the efficiency of Chebyshev-Legendre spectral method. The conclusion is given in section 6.

## 2. Basic Definitions and Fractional Derivatives

For $m \in \mathbb{N}$, the smallest integer that is greater than or equal to $\alpha$, i.e. $m = \lceil \alpha \rceil$, the Caputo's fractional derivative operator of order $\alpha > 0$, is defined as:

$$D^\alpha y(x) = \begin{cases} J^{m-\alpha} D^m y(x), & m-1 < \alpha \leq m, \\ D^m y(x), & \alpha = m, \end{cases} \tag{3}$$

where

$$J^{m-\alpha} y(x) = \frac{1}{\Gamma(m-\alpha)} \int_0^x (x-t)^{m-\alpha-1} y(t) dt, \quad v > 0, \quad x > 0.$$

For the Caputo's derivative we have [2]:

$$D^\alpha x^\beta = \begin{cases} 0, & \beta \in \{0,1,2,\ldots\} \text{ and } \beta < m, \\ \dfrac{\Gamma(\beta+1)}{\Gamma(\beta-\alpha+1)} x^{\beta-\alpha} & \beta \in \{0,1,2,\ldots\} \text{ and } \beta \geq m. \end{cases} \tag{4}$$

Recall that for $\alpha \in \mathbb{N}$, the Caputo differential operator coincides with the usual differential operator. Similar to standard differentiation, Caputo's fractional differentiation is a linear operator, i.e.,

$$D^\alpha (\lambda g(x) + \mu h(x)) = \lambda D^\alpha g(x) + \mu D^\alpha h(x),$$

where $\lambda$ and $\mu$ are constants.

The Chebyshev polynomials $\{T_i(t); i = 0,1,\ldots\}$ are defined on the interval $[-1,1]$ with the following recurrence formula:

$$T_{i+1}(t) = 2t\, T_i(t) - T_{i-1}(t), \quad i = 1,2,\ldots,$$

with $T_0(t) = 1$ and $T_1(t) = t$. The shifted Chebyshev polynomials are defined by introducing the change of variable $t = 2x - 1$. Let the shifted Chebyshev polynomials $T_i(2x - 1)$ be denote by $T_{1,i}(x)$, satisfying the relation

$$T_{1,i+1}(x) = 2(2x - 1)T_{1,i}(x) - T_{1,i-1}(x), \quad i = 1,2,\ldots, x \in [0,1], \tag{5}$$

where $T_{1,0}(x) = 1$ and $T_{1,1}(x) = 2x - 1$.

By these definitions we will have [32]

- $T_{1,i}(x) = i \sum_{k=0}^{i}(-1)^{i-k}\frac{(i+k-1)!\, 2^{2k}}{(i-k)!(2k)!} x^k \quad i = 1,2,\ldots.$ (6)
- $T_{1,i}(0) = (-1)^i, \quad T_{1,i}(1) = 1.$ (7)
- $\int_0^1 T_{1,j}(x) T_{1,k}(x)(x - x^2)^{\frac{-1}{2}} dx = \delta_{jk} h_k,$ (8)

where $\delta_{jk}$ is Kronecker delta and

$$h_k = \begin{cases} \pi, & k = 0 \\ \dfrac{\pi}{2}, & k \geq 1 \end{cases}, \tag{9}$$

In this paper, we will consider the Gauss-type quadrature formulas. We start by defining the Chebyshev-Gauss quadrature nodes and weights, respectively:

$$x_j = -\cos\frac{(2j+1)\pi}{2N+2}, \quad w_j = \frac{\pi}{N+1}, \quad j = 0,1,\ldots,N.$$

With the above choices, there holds

$$\int_{-1}^{1} p(x) \frac{1}{\sqrt{1 - x^2}} dx = \sum_{j=0}^{N} p(x_j) w_j, \quad \forall p \in P_{2N+1}, \tag{10}$$

where $P_{2N+1}$ is a polynomial of degree less than or equal $2N + 1$.

We now turn to the discrete Chebyshev transforms. The transforms can be performed via a matrix-vector multiplication with $O(N^2)$ operations as usual and when we use Chebyshev polynomials, it can be carried out with $O(N \log_2 N)$ operations via fast Fourier transform (**FFT**) [26-27].

We define the Chebyshev-Lagrange polynomial by

$$G_k(x) = \frac{T_k(x)}{(x - x_k)T_k'(x_k)}, \quad ;= 0,1,\ldots,N.$$

Given $u(x) \in C[-1,1]$, the Chebyshev-Lagrange interpolation operator $I_N^c u$ is defined by

$$(I_N^c u)(x) = \sum_{k=0}^{N} u_k G_k(x) \in \mathbb{P}_N, \tag{11}$$

where $\{u_k\}$ are determined by the forward discrete Chebyshev transform as follows

$$u_k = \sum_{j=0}^{N-1} u(x_j) \cos\frac{(2k+1)j\pi}{2N}, \quad 0 \leq k \leq N. \tag{12}$$

The above transform can be computed by using **FFT** in $O(N \log_2 N)$ operations [26-27].

Let $L_i(t)$ be the standard Legendre polynomial of degree i, then we have [20]

- Three-term recurrence relation
$$(i + 1)L_{i+1}(t) = (2i + 1)t\, L_i(t) - iL_{i-1}(t), \quad i \geq 1, \tag{13}$$
and the first two Legendre polynomials are

$$L_0(t) = 1, \quad L_1(t) = t.$$

- The Legendre polynomial $L_i(t)$ has the expansion

$$L_i(t) = \frac{1}{2^i} \sum_{l=0}^{[\frac{i}{2}]} (-1)^l \frac{(2i-2l)!}{2^l l!(i-l)!(i-2l)!} t^{i-2l}. \tag{14}$$

- Orthogonality

$$\int_{-1}^{1} L_j(t) L_k(t) \, dt = h_k \delta_{jk}, \tag{15}$$

such that

$$h_k = \frac{2}{2k+1}.$$

- Symmetry property

$$L_i(-t) = (-1)^i L_i(t), \quad L_i(\pm 1) = (\pm 1)^i. \tag{16}$$

Hence, $L_i(t)$ is an odd (resp. even) function, if $i$ is odd (resp. even).

Now, if we define the shifted Legendre polynomial of degree $i$ by $L_{1,i}(x) = L_i(2x-1)$, then we can obtain the analytic form and three-term recurrence relation of the shifted Legendre polynomials of degree $i$ by the following form, respectively

$$L_{1,i}(x) = \sum_{k=0}^{i} (-1)^{i+k} \frac{(i+k)!}{(i-k)!(k!)^2} x^k,$$

$$L_{1,i}(x) = \frac{2i+1}{i+1} x \, L_{1,i}(x) - \frac{i}{i+1} L_{1,i-1}(x), \quad i \geq 1. \tag{17}$$

According to Eq. (15), the orthogonality relation of shifted Legendre polynomials is

$$\int_0^1 L_{1,j}(t) L_{1,k}(t) \, dt = h_k \delta_{jk}. \tag{18}$$

We denote $L^2_\omega(I)$ by the weighted $L^2$ Hilbert space with the scalar product

$$(u,v) = \int_0^1 u(x) \, v(x) \, \omega(x) dt, \quad \forall u, v \in L^2_\omega(I),$$

and the norm $\|u\|_{L^2_\omega} = (u,u)^{\frac{1}{2}}_\omega$, where $\omega(x) = 1$ in the Legendre case and $\omega(x) = (1-x^2)^{-\frac{1}{2}}$ in the Chebyshev case. We may drop the subscript $w$ when $\omega = 1$. Therefore, the corresponding norm is

$$\|u\|_{L^2} = (u,u)^{\frac{1}{2}}.$$

Let $H^m_\omega(I) = \left\{ u \in L^2_\omega(I) : \frac{d^i u}{dx^i} \in L^2_\omega(I), \, i = 0,1,\ldots m \right\}$ be the weighted Sobolev space with the norm and semi norm defined respectively

$$\|y\|^2_{H^m_\omega(I)} = \sum_{k=0}^{m} \|y^{(k)}\|^2_{L^2_\omega(I)},$$

and

$$|y|^2_{H^{m:N}_\omega(I)} = \sum_{k=\min(m:N)}^{N} \|y^{(k)}\|^2_{L^2_\omega(I)}.$$

$H^m(I)$ by its inner product is Hilbert space.

For a function $y(x) \in L^2[0,1]$, the shifted Legendre expansion is
$$y(x) = \sum_{j=0}^{\infty} a_j L_{1,j}(x),$$
where
$$a_j = \frac{1}{h_{1,j}} \int_0^1 y(x) L_{1,j}(x) \, dx, \quad j = 0,1,2,\ldots, \tag{19}$$
and
$$h_{1,j} = \frac{1}{2} h_j = \frac{1}{2j+1}.$$

Now, we describe the Legendre-Gauss integration in the interval $(0,1)$. We denote by $x_{N,j}$, $\omega_{N,j}, j = 0,\ldots,N$, respectively the nodes and weights of the standard integration on the interval $(-1,1)$. We suppose $x_{1,N,j}, \omega_{1,N,j}, j = 0,\ldots,L$, are nodes and weights of the Legendre-Gauss integration in the interval $(0,1)$. Then, we have

$$x_{1,N,j} = \frac{1}{2}(x_{N,j} + 1), \quad w_{1,N,j} = \frac{1}{2} w_{N,j} \quad j = 0,\ldots,N.$$

According to Eq. (29) for any $g \in \mathbb{P}_{2N+1}$, set of all polynomials of degree at most $2N+1$, we get

$$\int_0^1 g(x) dx = \frac{1}{2} \int_{-1}^1 g\left(\frac{1}{2}(x+1)\right) dx$$
$$= \frac{1}{2} \sum_{j=0}^N \omega_{N,j} \, g\left(\frac{1}{2}(x_{N,j}+1)\right)$$
$$= \sum_{j=0}^N \omega_{1,N,j} \, g(x_{1,N,j}). \tag{20}$$

In practice, a number of first shifted Legendre polynomials are considered. We let

$$\phi(x) = [L_{1,0}(x), L_{1,1}(x), \ldots, L_{1,N}(x)]^T,$$
$$S_N(I) = span\{L_{1,0}(x), L_{1,1}(x), \ldots, L_{1,N}(x)\}. \tag{21}$$

**Theorem 2.1 [25]** suppose $\phi(x)$ is defined in Eq.(21) and $\alpha > 0$; then the following relation holds:

$$D^\alpha \phi(x) \cong D^{(\alpha)} \phi(x), \tag{22}$$

where $D^{(\alpha)}$ is the $(N+1) \times (N+1)$ operational matrix of Caputo derivative which is given by:

$$D^{(\alpha)} = (d_{ij})_{0 \leq i,j \leq N} = \begin{bmatrix} 0 & 0 & 0 & \cdots & 0 \\ & & \vdots & & \\ S_\alpha(m,0) & S_\alpha(m,1) & S_\alpha(m,2) & \cdots & S_\alpha(m,N) \\ & & \vdots & & \\ S_\alpha(i,0) & S_\alpha(i,1) & S_\alpha(i,2) & \cdots & S_\alpha(i,N) \\ & & \vdots & & \\ S_\alpha(N,0) & S_\alpha(N,1) & S_\alpha(N,2) & \cdots & S_\alpha(N,N) \end{bmatrix}, \qquad (23)$$

where

$$S_\alpha(i,j) = \sum_{k=m}^{i} \frac{(-1)^{i+k}\,(2j+1)\,(i+k)!\,\Gamma(k-j-\alpha+1)}{L^\alpha(i-k)!\,k!\,\Gamma(k-\alpha+1)\,\Gamma(k+j-\alpha+1)}. \qquad (24)$$

Note that because of $D^\alpha L_{1,i}(t) = 0$, for $i = 0,1,\ldots,m-1$, the first $m$ rows are zero in $\boldsymbol{D}$.

## 3. Chebyshev-Legendre Spectral Method

The Chebyshev-Legendre spectral method was introduced in [24] to take advantage of both the Legendre and Chebyshev polynomials. The main idea is to use the Legendre-Galerkin formulation which preserves the symmetry of the underlying problem and lead to a simple sparse linear system, while the physical values are evaluated at the Chebyshev-Gauss-type points. Thus, we may replace the expensive Legendre transform by a fast Chebyshev-Legendre transform between the coefficients of Legendre expansion and Chebyshev expansion at the Chebyshev-Gauss-type points.

The main advantage of using Chebyshev polynomials is that the discrete Chebyshev transform can be performed in $O(N \log_2 N)$ operations by using $FFT$. On the other hand, the discrete Legendre transform is expensive, and therefore in our article, the Chebyshev-Legendre method based on Legendre expansion and Chebyshev-Gauss-type points is applied to reduce the cost of solving the corresponding system (For more detail see [17, 24, 28]). Then, we use the Chebyshev interpolation operator $I_N^c$, relative to the Gauss-Chebyshev points to approximate the known functions and use of Legendre polynomials expansion to approximate the unknown function together. At last, the solution procedure is essentially the same as Legendre spectral method except that Chebyshev-Legendre transform, between the values of a function at the Gauss-Chebyshev points and the coefficients of its Legendre expansion, are needed instead of the Legendre transform. There are several efficient algorithms to transform from the coefficients of Legendre expansions to Chebyshev expansions at the Chebyshev-Gauss-Lobatto points and vice versa [24, 26-28]. We use the algorithm in [24] as follow:

We let

$$u(x) = \sum_{j=0}^{N} \alpha_j\, T_{1,j} = \sum_{j=0}^{N} \beta_j\, L_{1,j},$$
$$\boldsymbol{\alpha} = (\alpha_0, \alpha_1, \ldots, \alpha_N),$$
$$\boldsymbol{\beta} = (\beta_0, \beta_1, \ldots, \beta_N).$$

In this work, what we need to apply spectral method is using the transform between $\alpha$ and $\beta$. By virtue of orthogonality of Chebyshev and Legendre polynomials, the relation between $\alpha$ and $\beta$ can be obtained by computing $(u, T_{1,j})_w$ and $(u, L_{1,j})$. defining

$$A = (a_{ij})_{0 \le i,j \le N},$$
$$B = (b_{ij})_{0 \le i,j \le N},$$

then, by using Eqs. (8) and (15), we can obtain

$$a_{ij} = \frac{1}{h_i} (T_{1,i}, L_{1,j})_w,$$
$$b_{ij} = \left(i + \frac{1}{2}\right)(L_{1,i}, T_{1,j}).$$

Thus, we will have

$$\alpha = A\beta,$$
$$\beta = B\alpha,$$
$$AB = BA = I.$$

According to orthogonality and parity of the Chebyshev and Legendre polynomials, we get

$$a_{ij} = b_{ij} = 0, \text{ for } i > j \text{ or } i + j \text{ odd}.$$

Therefore, we only determine the nonzero elements of both $A$ and $B$ by using three-term recurrence relation of the shifted Legendre and Chebyshev polynomials. Applying definition of $a_{ij}$, we can obtain recurrence formula

$$\begin{aligned}
a_{ij} &= \frac{1}{h_i} (T_{1,i}, L_{1,j})_w \\
&= \frac{1}{h_i} \left(T_{1,i}, \frac{2j+1}{j+1}(2x-1)L_{1,j}(x) - \frac{j}{j+1}L_{1,j-1}(x)\right)_w \\
&= \frac{1}{h_i} \left\{\frac{2j+1}{j+1}\left((2x-1)T_{1,i}, L_{1,j}\right)_w - \frac{j}{j+1}(T_{1,i}, L_{1,j-1})_w\right\} \\
&= \frac{1}{h_i} \left\{\frac{2j+1}{2j+2}(T_{1,i+1} + T_{1,i-1}, L_{1,j})_w - \frac{j}{j+1} h_i a_{i,j-1}\right\} \\
&= \frac{h_{i+1}}{h_i} \frac{2j+1}{2j+2} a_{i+1,j} + \frac{h_{i-1}}{h_i} \frac{2j+1}{2j+2} a_{i-1,j} - \frac{j}{j+1} a_{i,j-1}.
\end{aligned}$$

We can similarly derive entries of matrix $B$ as follow

$$b_{ij} = \left(i + \frac{1}{2}\right) \tilde{b}_{ij},$$

where

$$\tilde{b}_{ij} := (L_{1,i}, T_{1,j}) = \frac{2i+2}{2i+1} \tilde{b}_{i+1,j} + \frac{2i}{2i+1} \tilde{b}_{i-1,j} - \tilde{b}_{i,j-1}.$$

Thus, we can obtain each nonzero element of $A$ and $B$ by just a few operations. Therefore, we can extremely apply Chebyshev-Legendre spectral method.

We now describe our spectral approximations to Eq. (1). Therefore, if $y_N(t) \in S_N(I)$, then by implementing Chebyshev-Legendre spectral method for Eq.(1), we can easily obtain

$$\sum_{i=0}^{n} a_i \left(D^{(i)} y_N, L_{1,k}\right) = (I_N^c f, L_{1,k}) + \left(\int_0^1 I_N^c k(.,s) \, D^\alpha y_N(s) \, ds, L_{1,k}\right). \tag{25}$$

We have $y_N(t) = \sum_{j=0}^{N} c_j L_{1,j}(t)$, then according to linearity of Caputo's fractional differentiation, Eq.(23) can be written as:

$$\sum_{i=0}^{n} a_i \sum_{j=0}^{N} c_j \left(D^{(i)} L_{1,j}, L_{1,k}\right)$$

$$= (I_N^c f, L_{1,k}) + \sum_{j=0}^{N} c_j \left(\int_0^1 I_N^c k(.,s) \, D^\alpha L_{1,j}(s) \, ds, L_{1,k}\right). \tag{26}$$

From Eq. (22) to (24) in **Theorem 2.1**, we can obtain

$$D^\alpha L_{1,j}(t) = \sum_{l=0}^{N} S_\alpha(j,l) L_{1,l}(t), \qquad j = m, m+1, \dots, N. \tag{27}$$

We notice that if $\alpha = n \in \mathbb{N}$, then $S_\alpha$ defined in Eq. (24) tend to integer order case and **Theorem 2.1** gives the same result as integer order case.

Inserting Eq.(27) in Eq.(26), we get

$$\sum_{i=0}^{n} \sum_{j=i}^{N} a_i c_j \left(\sum_{l=0}^{N} S_i(j,l) (L_{1,l}, L_{1,k})\right)$$

$$= (I_N^c f, L_{1,k}) + \sum_{j=m}^{N} c_j \sum_{l=0}^{N} S_\alpha(j,l) \left(\int_0^1 I_N^c k(.,s) \, L_{1,l}(t) \, ds, L_{1,k}\right). \tag{28}$$

Then, making use of the orthogonality relation of shifted Legendre polynomials, i.e. Eq.(18), Eq. (24) reduce to

$$\sum_{i=0}^{n} \sum_{j=i}^{N} a_i\, c_j \frac{S_i(j,k)}{2k+1}$$

$$= (I_N^c f, L_{1,k}) + \sum_{j=m}^{N} c_j \sum_{l=0}^{N} S_\alpha(j,l) \left( \int_0^1 I_N^c k(.,s)\, L_{1,l}(t)\, ds\, , L_{1,k} \right). \tag{29}$$

We let

$$h_l(x) = \int_0^1 I_N^c k(x,s)\, L_{1,l}(t)\, ds \cong \sum_{r=0}^{N} b_{lr}\, L_{1,r}(x),$$

$$f_k = (I_N^c f, L_{1,k}).$$

Thus, again by using Eqs.(18), Eq. (29) becomes the following form

$$\sum_{i=0}^{n} \sum_{j=i}^{N} a_i\, c_j \frac{S_i(j,k)}{2k+1} = f_k + \sum_{j=m}^{N} \sum_{l=0}^{N} c_j\, S_\alpha(j,l) \frac{b_{lk}}{2k+1}. \tag{30}$$

It is easy to verify that initial conditions convert to following equations

$$\sum_{j=0}^{N} \sum_{l=0}^{N} c_j\, S_i(j,l)\, L_{1,l}(0) = d_i, \quad i = 0,1,\ldots,n-1. \tag{31}$$

Combining Eqs. (30) and (31) yields

$$\begin{cases} \displaystyle\sum_{i=0}^{n} \sum_{j=i}^{N} a_i\, c_j \frac{S_i(j,k)}{2k+1} - \sum_{l=0}^{N} \sum_{j=m}^{N} c_j\, S_\alpha(j,l) \frac{b_{lk}}{2k+1} = f_k, & k = 0,1,\ldots N-n, \\ \displaystyle\sum_{j=0}^{N} \sum_{l=0}^{N} c_j\, S_i(j,l)\, L_{1,l}(0) = d_i, & i = 0,1,\ldots, n-1. \end{cases}$$

By solving the above system of linear equations, we can get the value of $\{c_j\}_{j=0}^{N}$ and obtain the expression of $y_N(x)$ accordingly.

## 4. Convergence Analysis of the Chebyshev-Legendre Spectral method

In this section, we present a general approach to the convergence analysis for NIFDEs that is proved in $L^2-$norm. Here, there are some properties and elementary lemmas, which are important for the derivation of the main results.

**Lemma 4.1 [29]** For multiple integrals, the following relation holds:

$$\int_0^t \int_0^{t_n} \ldots \int_0^{t_3} \int_0^{t_2} g(t_1)\, dt_1 dt_2 \ldots dt_n = \frac{1}{(n-1)!} \int_0^t (t-s)^{n-1} g(s)\, ds, \tag{32}$$

where $g$ is integrable function on interval $(0,t)$ and $t_i$ $(i=2,3,\ldots,n)$ are parameters in the purpose interval.

**Lemma 4.2 [30]** (Granwall's Lemma) Assume that $u, \omega, \beta \in C(I)$ with $\beta(t) \geq 0$. If $u$ satisfies the inequality

$$u(t) \leq \omega(t) + \int_0^t \beta(s) u(s) \, ds, \quad t \in I,$$

then

$$u(t) \leq \omega(t) + \int_0^t \beta(s) \omega(s) \exp\left(\int_s^t \beta(v) \, dv\right) ds, \quad t \in I. \tag{33}$$

On the other word, if $\omega$ is non-decreasing on $I$, the above inequality reduce to

$$u(t) \leq \omega(t) \exp\left(\int_s^t \beta(v) \, dv\right), \quad t \in I. \tag{34}$$

**Lemma 4.3 [30]** Suppose that $k$ is a given kernel function on $I \times I$. If $f \in L^p(a,b)$ for $1 \leq p \leq \infty$, the integral

$$Tf(x) = \int_a^{x \text{ or } b} k(x,t) f(t) \, dt$$

is well-defined in $L^p(a,b)$ and there exists $C^* > 0$ such that

$$\|Tf\|_{L^p(a,b)} \leq C^* \|f\|_{L^p(a,b)}. \tag{35}$$

Let $p_N$ be the interpolation projection operator from $\mathbb{L}^2(I)$ upon $\mathbb{P}_N(I)$. Then, for any function $f$ in $L^2(I)$ satisfies

$$\int_0^1 (f - p_N f)(t) \, g(t) dt = 0, \quad \forall g \in \mathbb{P}_N(I).$$

Also, the following relations for interpolation in shifted Legendre polynomials and shifted Gauss-Legendre nodal points $k \geq 1$ (or for any fixed $k \leq N$) may readily be obtained as [14]

$$\|y - p_N(y)\|_{H^l(I)} \leq C_1 N^{2l - \frac{1}{2} - k} |y|_{H^{k:N}(I)}, \tag{36}$$

$$\|y - p_N(y)\|_{L^2(I)} \leq C_2 N^{-k} |y|_{H^{k:N}(I)}. \tag{37}$$

where $y \in H^m(I)$, and $C_1$ and $C_2$ are constants independent of $N$ and $0 \leq l \leq m$.

Now, we shall prove the main result in this section. In the following theorem, an error estimation for an approximate solution of Eq. (1) with supplementary conditions of Eq. (2) is obtained. Let $e_N(x) = y(x) - y_N(x)$, be the error function of the Chebyshev-Legendre spectral approximation to $y(x)$. From the mathematical point of view, it is possible to keep track of the effect of the boundary conditions upon the overall accuracy of the scheme. In the other hand, the boundary treatment does not destroy the spectral accuracy of the Chebyshev-Legendre method.

**Theorem 4.3** For sufficiently large $N$, the Chebyshev-Legendre spectral approximations $y_N(x)$ converge to exact solution in $L^2$-norm, i.e.

$$\|e_N\|_{L^2(I)} = \|y - y_N\|_{L^2(I)} \to 0.$$

Proof. Assume that $y_N(x)$ is obtained by using the Chebyshev-Legendre spectral method Eq. (1) together with initial conditions Eq. (2). Then, we have

$$\sum_{i=0}^{n} a_i D^{(i)} y_N(t) = p_N(f(t)) + p_N\left(\int_0^1 k(.,s) D^\alpha y_N(s) \, ds\right), \tag{38}$$

such that $p_N$ is the Lagrange interpolation polynomial operator defined for Legendre polynomial. With $n$ times integration from Eq. (38), we obtain

$$\sum_{i=0}^{n} a_i \int_0^t \int_0^{t_n} \cdots \int_0^{t_3} \int_0^{t_2} y_N^{(i)}(t_1) \, dt_1 dt_2 \ldots dt_n$$

$$= \int_0^t \int_0^{t_n} \cdots \int_0^{t_3} \int_0^{t_2} p_N(f(t_1)) \, dt_1 dt_2 \ldots dt_n$$

$$+ \int_0^t \int_0^{t_n} \cdots \int_0^{t_3} \int_0^{t_2} p_N\left(\int_0^1 k(t_1,s) D^\alpha y_N(s) \, ds\right) dt_1 dt_2 \ldots dt_n. \tag{39}$$

By virtue of **Lemma 4.1,** we can convert each of the multiple integral to single integral, so we have

$$a_n y_N(t) + g(t) + \sum_{i=0}^{n-1} \int_0^t \frac{a_i}{(n-i-1)!} (t-s)^{n-i-1} y_N(s) \, ds$$

$$= \int_0^t \frac{(t-s)^{n-1}}{(n-1)!} p_N(f(s)) \, ds$$

$$+ \int_0^t \frac{(t-s)^{n-1}}{(n-1)!} p_N\left(\int_0^1 k(s,s_1) D^\alpha y_N(s_1) \, ds_1\right) ds, \tag{40}$$

where $g$ is a polynomial of degree $n$ with the initial condition coefficient. Similarly, from Eq. (1), we get

$$a_n y(t) + g(t) + \sum_{i=0}^{n-1} \int_0^t \frac{a_i}{(n-i-1)!} (t-s)^{n-i-1} y(s) \, ds$$

$$= \int_0^t \frac{(t-s)^{n-1}}{(n-1)!} f(s) \, ds$$

$$+ \int_0^t \frac{(t-s)^{n-1}}{(n-1)!} \int_0^1 k(s,s_1) D^\alpha y(s_1) \, ds_1 \, ds. \tag{41}$$

By subtracting Eq. (40) from Eq. (41), we obtain

$$a_n e_N(t) + \sum_{i=0}^{n-1} \int_0^t \frac{a_i}{(n-i-1)!} (t-s)^{n-i-1} e_N(s)\, ds$$
$$= \int_0^t \frac{(t-s)^{n-1}}{(n-1)!} e_{p_N}(f(s))\, ds + \int_0^t \frac{(t-s)^{n-1}}{(n-1)!} e_{p_N}(K_\alpha y(s))\, ds, \qquad (42)$$

such that

$$e_{p_N}(f(s)) = f(s) - p_N(f(s)),$$
$$K_\alpha y(s) = \int_0^1 k(s, s_1)\, D^\alpha y(s_1)\, ds_1,$$
$$e_{p_N}(K_\alpha y(s)) = K_\alpha y(s) - p_N(K_\alpha y_N(s))$$
$$= K_\alpha y(s) - K_\alpha y_N(s) + K_\alpha y_N(s) - p_N(K_\alpha y_N(s))$$
$$= K_\alpha e_N(s) - e_{p_N}(K_\alpha y_N(s)). \qquad 0 \le s \le t.$$

From Eq. (42), we can obtain

$$|e_N(t)| \le \sum_{i=0}^{n-1} \left| \frac{a_i}{-a_n\, (n-i-1)!} \right| \left| \int_0^t |(t-s)^{n-i-1} e_N(s)| ds \right.$$
$$+ \frac{1}{|a_n|(n-1)!} \int_0^t |(t-s)^{n-1} e_{p_N}(f(s))|\, ds$$
$$+ \frac{1}{|a_n|(n-1)!} \int_0^t |(t-s)^{n-1} e_{p_N}(K_\alpha y(s))|\, ds$$
$$\le C_2 \int_0^t |e_N(s)|\, ds + C_3 \int_0^t |e_{p_N}(f(s))|\, ds + C_4 \int_0^t |e_{p_N}(K_\alpha y(s))|\, ds. \qquad (43)$$

Applying **Lemma 4.2** leads to

$$|e_N(t)| \le \exp\left( \int_0^t C_2\, ds \right) \left( C_3 \int_0^t |e_{p_N}(f(s))|\, ds + C_4 \int_0^t |e_{p_N}(K_\alpha y(s))|\, ds \right)$$
$$\le C_5 \int_0^t |e_{p_N}(f(s))|\, ds + C_6 \int_0^t |e_{p_N}(K_\alpha y(s))|\, ds. \qquad (44)$$

Equivalently, by using the $L^2$-norm, we get

$$\|e_N\|_{L^2(I)} \le C_5 \left\| \int_0^t |e_{p_N}(f(s))|\, ds \right\|_{L^2(I)} + C_6 \left\| \int_0^t |e_{p_N}(K_\alpha y(s))|\, ds \right\|_{L^2(I)}. \qquad (45)$$

Bu using **Lemma 4.3,** the above inequality reduce to

$$\|e_N\|_{L^2(I)} \leq C_7 \left\|e_{p_N}(f(s))\right\|_{L^2(I)} + C_8 \left\|e_{p_N}(K_\alpha y(s))\right\|_{L^2(I)}. \tag{46}$$

On the other hand, we have

$$\begin{aligned}
e_{p_N}(K_\alpha y(s)) &= \int_0^1 k(s,s_1) D^\alpha y(s_1) \, ds_1 - p_N \left( \int_0^1 k(s,s_1) D^\alpha y_N(s_1) \, ds_1 \right) \\
&= \int_0^1 k(s,s_1) D^\alpha y(s_1) \, ds_1 - \int_0^1 k(s,s_1) D^\alpha y_N(s_1) \, ds_1 \\
&\quad + \int_0^1 k(s,s_1) D^\alpha y_N(s_1) \, ds_1 - p_N \left( \int_0^1 k(s,s_1) D^\alpha y_N(s_1) \, ds_1 \right) \\
&= \int_0^1 k(s,s_1) D^\alpha e_N(s_1) \, ds_1 + E(s),
\end{aligned} \tag{47}$$

where

$$E(s) = \int_0^1 k(s,s_1) D^\alpha y_N(s_1) \, ds_1 - p_N \left( \int_0^1 k(s,s_1) D^\alpha y_N(s_1) \, ds_1 \right).$$

Therefor

$$\left\|e_{p_N}(K_\alpha y(s))\right\|_{L^2(I)} \leq \left\| \int_0^1 k(.,s_1) D^\alpha e_N(s_1) \, ds_1 \right\|_{L^2(I)} + \|E(s)\|_{L^2(I)}. \tag{48}$$

According to relation (37) and **Lemma 4.3,** we have

$$\begin{aligned}
\|E(s)\|_{L^2(I)} &\leq C_2 N^{-k} \left\| \int_0^1 k(.,s_1) D^\alpha y_N(s_1) \, ds_1 \right\|_{H^{k:N}(I)} \\
&\leq C_2 C^* N^{-k} \|D^\alpha y\|_{H^{k:N}(I)}
\end{aligned} \tag{49}$$

In the other hand, because of linear operator $D^\alpha : \mathbb{P}_N \to \mathbb{P}_N$ is continuous and bounded [31], thus there exists a constant $C^{**} \geq 0$ such that

$$\|D^\alpha y_N\|_{H^{k:N}(I)} \leq C^{**} \|y_N\|_{H^{k:N}(I)}, \quad k \in \mathbb{N} \text{ and } k \leq N. \tag{50}$$

Therefore, by combining two recent relations, we get

$$\begin{aligned}
\|E(s)\|_{L^2(I)} &\leq C_2 C^* C^{**} N^{-k} \|y_N\|_{H^{k:N}(I)} \\
&= C_9 N^{-k} \|y - e_N\|_{H^{k:N}(I)}
\end{aligned}$$

$$\leq C_9 \, N^{-k} \left( \|y\|_{H^{k:N}(I)} + \|e_N\|_{H^{k:N}(I)} \right). \tag{51}$$

Now, by using relation (36), we proceed with the above inequality as

$$\begin{aligned}
\|E(s)\|_{L^2(I)} &\leq C_9 \, N^{-k} \left( \|y\|_{H^{k:N}(I)} + \|e_N\|_{H^{k:N}(I)} \right) \\
&\leq C_9 \, N^{-k} \left( \|y\|_{H^{k:N}(I)} + \|e_N\|_{H^{1:N}(I)} \right) \\
&\leq C_9 \, N^{-k} \left( \|y\|_{H^{k:N}(I)} + C_1 N^{\frac{3}{2}-k} |y|_{H^{k:N}(I)} \right) \\
&\leq C_9 \, N^{-k} \|y\|_{H^{k:N}(I)} + C_{10} \, N^{\frac{3}{2}-2k} |y|_{H^{k:N}(I)}. 
\end{aligned} \tag{52}$$

Similarly, from **Lemma 4.3** and relation (50), we obtain

$$\begin{aligned}
\left\| \int_0^1 k(\cdot,s_1) \, D^\alpha e_N(s_1) \, ds_1 \right\|_{L^2(I)} &\leq C_{11} \, \|D^\alpha e_N\|_{H^{k:N}(I)} \\
&\leq C_{11} \, \|e_N\|_{H^{1:N}(I)} \\
&\leq C_{12} \, N^{\frac{3}{2}-k} |y|_{H^{k:N}(I)}.
\end{aligned} \tag{53}$$

At last, combining (48), (52), and (53) gives

$$\left\| e_{p_N}(K_\alpha y(s)) \right\|_{L^2(I)} \leq C_9 \, N^{-k} \|y\|_{H^{k:N}(I)} + C_{13} \, N^{\frac{3}{2}-2k} |y|_{H^{k:N}(I)}. \tag{54}$$

In a similar manner with relation (37), we may write

$$\left\| e_{p_N}(f(s)) \right\|_{L^2(I)} \leq C_2 N^{-k} |f|_{H^{k:N}(I)}. \tag{55}$$

Finally, by substituting $(54) - (55)$ in (46), the following relation can be obtained

$$\begin{aligned}
\|e_N\|_{L^2(I)} &\leq C_7 \, C_2 N^{-k} |f|_{H^{k:N}(I)} + C_8 \left( C_9 \, N^{-k} \|y\|_{H^{k:N}(I)} + C_{13} \, N^{\frac{3}{2}-2k} |y|_{H^{k:N}(I)} \right) \\
&\leq \gamma_1 \, N^{-k} \left( |f|_{H^{k:N}(I)} + \|y\|_{H^{k:N}(I)} \right) + \gamma_2 N^{\frac{3}{2}-2k} |y|_{H^{k:N}(I)}.
\end{aligned}$$

The above inequality proves that the approximation is convergent in $L^2$-norm. Hence the theorem is proved.

## 5. Numerical example

To show efficiency of the numerical method, the following examples are considered.

**Example 5.1.** Consider the following fractional integro-differential equation [20]

$$y'(x) = 14\left(1 - \frac{t}{2.5\,\Gamma(1.5)}\right) + \int_0^1 xs\, D^{\frac{1}{2}}y(s)\, ds,$$

with the initial condition: $y(0) = 0$, and exact solution $y(x) = 14x$.

We have solved this example using Chebyshev-Legendre spectral method and approximations are obtained as follows:

$$\begin{aligned} n = 0: &\quad y_0(x) = 0, \\ n = 1: &\quad y_1(x) = 14x, \\ n = 2: &\quad y_2(x) = 14x, \\ n = 3: &\quad y_3(x) = 14x, \end{aligned}$$

and so on. Therefore, we obtain $y(x) = 14x$ which is the exact solution of the problem.

**Example 5.2.** Our second example is the following fractional integro-differential equation [20]

$$y'(x) = f(x) + \int_0^1 x^2 s^2\, D^{\frac{1}{4}}y(s)\, ds,$$

with the initial condition

$$y(0) = 0,$$

where $f(x) = 8x^3 - \frac{3}{2}x^{\frac{1}{2}} - \left(\frac{48}{6.75\,\Gamma(4.75)} - \frac{\Gamma(2.75)}{4.25\,\Gamma(2.25)}\right)x^2$, and $y(x) = 2x^4 - x^{\frac{3}{2}}$ is the exact solution. The numerical results of our method can be seen from **Figure 1** and **Figure 2**. These results indicate that the spectral accuracy is obtained for this problem, although the given function $f(t)$ is not very smooth.

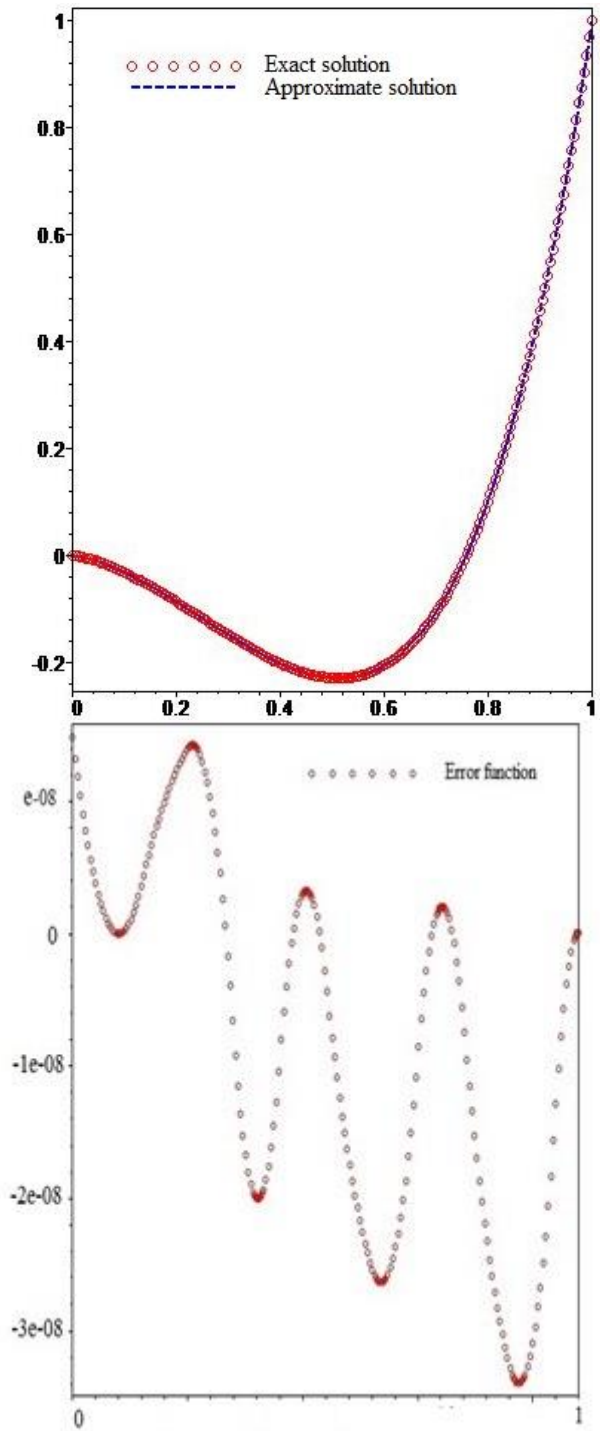

**Figure 1.** Comparison between exact solution and approximate solution of Example 5.2 (left), the error function for some different values (right)

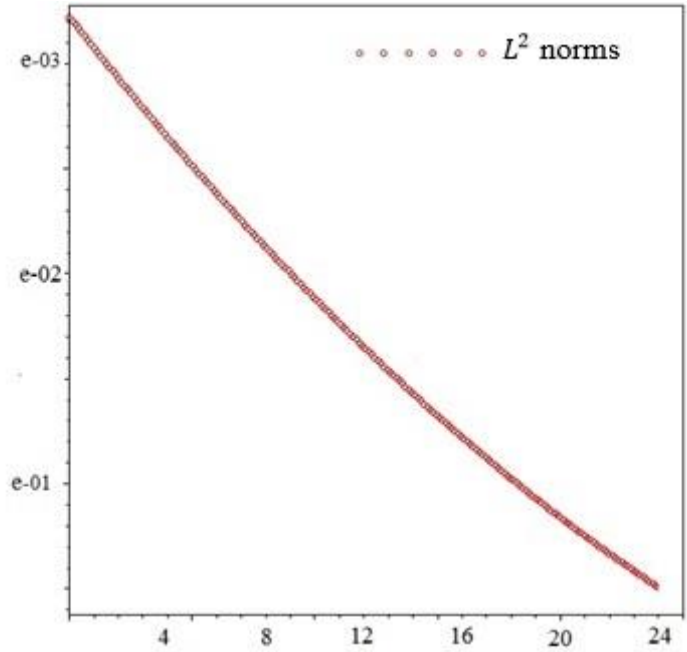

**Figure 2.** The error of numerical and exact solution of Example 5.2 versus the number of interpolation operator in $L^2$ norm

**Example 5.3.** Consider the following fractional integro-differential equation

$$2y''(x) + y'(x) = \left(\frac{9\sqrt{\pi} - 12}{\sqrt{\pi}}\right)x^2 + 36x + 8 + \int_0^1 x^2\sqrt{s}\, D^{\frac{3}{2}}y(s)\, ds,$$

with the initial conditions

$$y(0) = 0,$$
$$y'(0) = 8.$$

Taking $N = 4$, by implementing the Chebyshev-Legendre spectral method, we get the numerical solution as follow

$$y_4 = 8x + 1.003417 \times 10^{-12} x^2 + 3x^3 + 1.652105 \times 10^{-13} x^4 \cong 8x + 3x^3.$$

The approximate solution $y_4$ for this fractional integro-differential equation tends rapidly to exact solution, i.e. $y(x) = 8x + 3x^3$.

**Example 5.4.** Let us consider the following fractional integro-differential equation

$$3y^{(3)}(x) - y''(x) + y(x) = \left(7 - \frac{32}{15\sqrt{\pi}}\right)e^x + 3xe^x + \int_0^1 e^{x-s} D^{\frac{1}{2}}y(s)\, ds,$$

with the initial conditions

$$y(0) = 0,$$
$$y'(0) = 1,$$

$$y''(0) = 2.$$

The exact solution of this fractional integro-differential is $y(x) = xe^x$. We have reported the obtained numerical results for $N = 4$ and 8 in **Table 1**. Also, in **Figure 3**, we plot the resulting errors versus the number $N$ of the steps. This figure shows the exponential rate of convergence predicted by the proposed method.

| t | Proposed method at $N = 4$ | Proposed method at $N = 8$ | Exact solution |
|---|---|---|---|
| 0 | .0000000000 | .0000000000 | .0000000000 |
| 0.2 | .2442815491 | .2442805512 | .2442805516 |
| 0.4 | .5967277463 | .5967298754 | .5967298792 |
| 0.6 | 1.093273679 | 1.093271221 | 1.093271280 |
| 0.8 | 1.780432013 | 1.780432791 | 1.780432742 |
| 1 | 2.718281658 | 2.718281815 | 2.718281828 |

**Table 1.** The numerical results and exact solution for Example 5.4

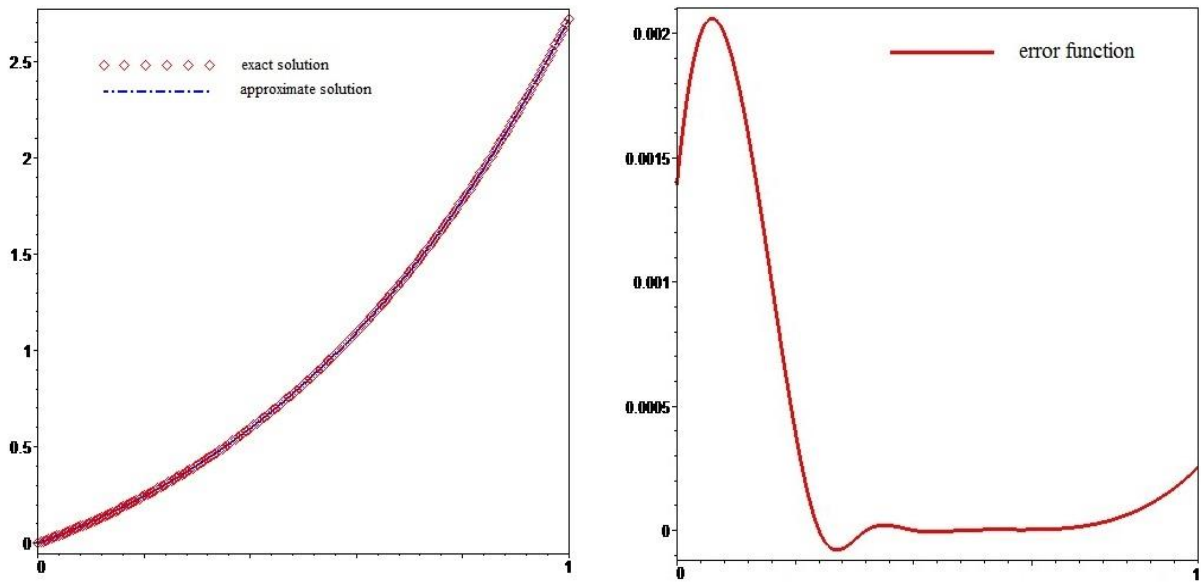

**Figure 3.** Comparison between exact solution and approximate solution of Example 5.4 (left), the error function for some different values (right)

## 6. Conclusion

In this paper we present a Chebyshev-Legendre spectral approximation of a class of Fredholm fractional integro-differential equations. The most important contribution of this paper is that the errors of approximations decay exponentially in $L^2 - norm$. We prove that our proposed method is effective and has high convergence rate. The results given in the previous section are compared with exact solutions. The satisfactory results agree very well with exact solutions only for small numbers of shifted Legendre polynomials.

**Reference**

[1] A.A. Kilbas, H.M. Srivastava and J.J. Trujillo, Theory and Applications of Fractional Differential Equations, North-Holland Math. Stud., vol. 204, Elsevier Science B.V., Amsterdam, (2006).

[2] I. Podlubny, Fractional Differential Equations, Academic Press, San Diego, (1999).

[3] A. H. Bhrawy and M. A. Alghamdi, A shifted Jacobi-Gauss-Lobatto collocation method for solving nonlinear fractional Langevin equation involving two fractional orders in different intervals, Boundary Value Problems, vol. 2012, article 62, 13 pages, 2012.

[4] E. Rawashdeh, Numerical solution of fractional integro-differential equations by collocation method, Appl. Math. Comput. 176, 1-6 (2006).

[5] Y. Nawaz, Variational iteration method and Homotopy perturbation method for fourth-order fractional integro-differential equations, Comput. Math. Appl. 61, 2330-2341 (2011).

[6] S. Momani and M. A. Noor, Numerical method for fourth-order fractional integro-differential equations, Appl. Math. Comput. 182, 754-760 (2006).

[7] H. Saeedi and F. Samimi, He's homotopy perturbation method for nonlinear Fredholm integro-differential equation of fractional order, International J. of Eng. Res. And App., Vol. 2, no. 5, pp. 52-56, 2012.

[8] A. Ahmed and S. A. H. Salh, Generalized Taylor matrix method for solving linear integro-fractional differential equations of Volterra type, App. Math. Sci., Vol. 5, no. 33-36, pp. 1765-1780, 2011.

[9] A. Arikoglo and L. Ozkol, Solution of fractional integro-differential equation by using fractional differential transform method, Chaos Solitons Fractals 34, 1473-1481 (2007).

[10] J. Wei and T. Tian, Numerical solution of nonlinear Volterra integro-differential equations of fractional order by the reproducing kernel method, Appl. Math. Model. 39, 4871-4876 (2015).

[11] Li Zhu and Qibin Fan, Solving fractional nonlinear Fredholm integro-differential equations by the second kind Chebyshev wavelet, Common. Nonlinear Sci. Numer. Simul. 17, 2333–2341(2012).

[12] H. Saeedi, M. Mohseni Moghadam, N. Mollahasani and GN. Chuev, A CAS wavelet method for solving nonlinear Fredholm integro-differential equation of fractional order, Common. Nonlinear Sci. Numer. Simul. 16, 1154-1163 (2011).

[13] Z. Meng, L. Wang, H. Li and W. Zhang, Legendre wavelets method for solving fractional integro-differential equations, Int. J. Comput. Math. 92, 1275-1291 (2015).

[14] C. Canuto, M. Y. Hussaini, A. Quarteroni and T. A. Zang, Spectral methods: Fundamentals in single domain, Springer, Berlin (2006).

[15] P. E. Bjorstad and B. P. Tjostheim, Efficient algorithm for solving a fourth-order equation with the spectral Galerkin method, SIAM J. Sci. Comput., 18, pp. 621-632 (1997).

[16] J. Shen, Efficient spectral-Galerkin method I. Direct solvers for second- and fourth-order equations by using Legendre polynomials, SIAM J. Sci. Comput., 15, pp. 1489-1505 (1994).

[17] J. Shen, T. Tang and Li-Lian Wang, Spectral Methods: Algorithms, Analysis and Applications, Springer-Verlag, Berlin Heidelberg (2011).

[18] S. Karimi Vanani and A. Amin Ataei, Operational Tau approximation for a general class of fractional integro-differential equations, J. Comp. and Appl. Math., Vol. 30, N. 3, 655-674 (2011).

[19] M. M. Khader and N. H. Sweilam, A Chebyshev pseudo-spectral method for solving fractional order integro-differential equation, The ANZIAM J. of Australian Math. Soc., Vol. 51 (04), pp. 464-475 (2010).

[20] A. Yousefi, T. Mahdavi-Rad and S. G. Shafiei, A quadrature Tau method for solving Fractional integro-differential equations in the Caputo sense, J. of Math. Computer Sci. 16, 97-107 (2015).

[21] W. G. El-Sayed and A. M.A. El-Sayed, On the fractional integral equations of mixed type integro-differential equations of fractional orders, App. Math. Comput. 154, 461-467 (2004).

[22] M. F. Al-Jamal and E. A. Rawashdeh, The approximate solution of fractional integro-differential equations, Int. J. Contemp. Math. Sci. 4, no. 22, 1067-1078 (2009).

[23] J. Shen and R. Temam, Nonlinear Galerkin method using Chebyshev and Legendre polynomials I. The one-dimensional case, SIAM J. Numer. Anal., 32(1995), pp. 215-234.

[24] J. Shen, Efficient Chebyshev-Legendre Galerkin methods for elliptic problems, in Proc. ICOSAHOM'95, A. V. Ilin and R. Scott, eds., Houston J. Math., 233-240 (1996).

[25] A. H. Bhrawy, A.S. Alofi and S.S. Ezz-Eldien, A quadrature Tau method for fractional differential equations with variable coefficients, App. Math. Letter 24, 2146-2152 (2011).

[26] B.K. Alpert and V. Rokhlin, A fast algorithm for the evaluation of Legendre expansions. SIAM J. Sci. Stat. Comput., 12:158-179 (1991).

[27] L. Greengard and V. Rokhlin. A fast algorithm for particle simulations. J. Comput. Phys., 73:325-348 (1987).

[28] W.S. Don and D. Gottlieb. The Chebyshev-Legendre method: implementing Legendre methods on Chebyshev points. SIAM J. Numer. Anal., 31:1519-1534 (1994).

[29] R.P. Kanwal, Linear Integral Equations, second ed., Birkhauser Boston, Inc., Boston, MA, (1971).

[30] Y. Chen and T. Tang, Convergence analysis of the Jacobi spectral collocation methods for Volterra integral equations with a weakly singular kernel, Math. Comp. 79 (269): 147–167(2010).

[31] A.W. Naylor and G.R. Sell, Linear operator theory in engineering and science, Second ed., in: Applied mathematical sciences, Vol. 40, Springer-Verlag, New York, Berlin (1982).

[32] E.H. Doha, A.H. Bhrawy and S.S. Ezz-Eldien, A Chebyshev spectral method based on operational matrix for initial and boundary value problems of fractional order, J. of Comp. and Math. with Appl. 62 2364–2373 (2011).